\documentclass[10pt]{amsart}
\theoremstyle{plain}
\usepackage{amsmath, amsfonts, epsfig}
\usepackage[all]{xy}
\title[Cyclic Branched Covers]{Combinatorial Description of Knot Floer Homology of Cyclic Branched Covers}
\author{J. Elisenda Grigsby} 
\thanks{The author was partially supported by an NSF Postdoctoral Fellowship.}
\address{Columbia Math Dept.;2990 Broadway MC4406; NY, NY 10027}
\email{egrigsby@math.columbia.edu}

\begin{document}
\bibliographystyle{halpha}

\newtheorem{theorem}{Theorem}[section]
\newtheorem{lemma}{Lemma}[section]
\newtheorem{definition}{Definition}[section]
\newtheorem{proposition}{Proposition}[section]
\newtheorem{corollary}{Corollary}[section]
\newcommand{\Ozsvath}{{Ozsv{\'a}th} }
\newcommand{\Szabo}{{Szab{\'o}} }

\begin{abstract}
In this paper, we introduce a simple combinatorial method for computing all versions ($\wedge,+,-,\infty$) of the knot Floer homology of the preimage of a two-bridge knot $K_{p,q}$ inside its double-branched cover, $-L(p,q)$.  The $4$-pointed genus $1$ Heegaard diagram we obtain looks like a twisted version of the toroidal grid diagrams recently introduced by Manolescu, Ozsv{\'a}th, and Sarkar.  We conclude with a discussion of how one might obtain nice Heegaard diagrams for cyclic branched covers of more general knots.
\end{abstract}

\maketitle
\section{Introduction}\label{section:Introduction}

Heegaard Floer homology, introduced by \Ozsvath and \Szabo in \cite{MR2113019}, associates to a closed, oriented, connected three-manifold $Y$ and a Spin$^c$ structure $\mathfrak{s} \in$ Spin$^c(Y)$ a collection, $$\widehat{HF}(Y;\mathfrak{s}), HF^\infty(Y;\mathfrak{s}), HF^-(Y;\mathfrak{s}), HF^+(Y;\mathfrak{s}),$$ of graded abelian groups, most naturally thought of as the homology groups of chain complexes with coefficients in $\mathbb{Z}$, $\mathbb{Z}[U,U^{-1}]$, $\mathbb{Z}[U]$, and $\frac{\mathbb{Z}[U,U^{-1}]}{U \cdot \mathbb{Z}[U,U^{-1}]}$, respectively.

The additional data of a nullhomologous, oriented link $L$ in $Y$ provides refinements of these invariants, discovered by \Ozsvath and \Szabo \cite{MR2065507} and independently by Rasmussen \cite{GT0306378}, that have had remarkable success in providing new information about links and the three-manifolds obtained by surgery upon them.  See, for example, \cite{MR2026543,MR2023281,MR2153455, MR2168576}.

The power of these invariants has proved difficult to harness, however, since their computation necessitates counts of holomorphic disks.  Recent work of Manolescu, {Ozsv{\'a}th}, Sarkar, and Wang \cite{GT0607691, GT0607777} has markedly improved the situation by providing a purely combinatorial definition of the invariant in most cases.  Specifically, \cite{GT0607691} provides a combinatorial description of all versions of knot Floer homology for a knot $K$ in $S^3$ while \cite{GT0607777} provides one for the filtered chain homotopy type of the complex $\widehat{CF}(Y)$ in the presence of a nullhomologous knot $K \subset Y$, $Y$ any (closed, oriented, connected) three-manifold.

Although algorithms now exist to perform previously inaccessible calculations, the chain complexes arising from the general algorithms are quite large.

The aim of the present paper is to provide a streamlined combinatorial description of the knot Floer homology groups of the preimage, $\widetilde{K}_{p,q}$, of a two-bridge knot $K_{p,q} \subset S^3$ inside its $2$-fold cyclic branched cover, $\Sigma^2(K)= -L(p,q)$.  

Interest in studying Heegaard Floer homology in cyclic branched covers of a knot $K$ has so far centered upon obtaining new concordance invariants \cite{GT0508065}, still a promising direction for further exploration.  More generally, it is shown in \cite{GT0507498} that $\widehat{HFK}(\Sigma^n(K);\widetilde{K})$ captures nonabelian information about $\pi_1(S^3-K)$ unavailable through the study of $\widehat{HFK}(S^3;K)$, suggesting the likelihood of further applications.

In brief, we obtain a $4$-pointed genus $1$ Heegaard diagram compatible with $\widetilde{K}_{p,q} \subset \Sigma^2(K_{p,q})$ which is a {\it twisted toroidal grid diagram} consisting of two parallel curves of slope $0$ and two of slope $\frac{p}{q}$, partitioning the torus into $2pq$ cells.

More specifically, we identify the universal cover of the torus with the plane: $$T^2 := \mathbb{R}^2/\mathbb{Z}^2.$$

\begin{figure}
\begin{center}
\epsfig{file=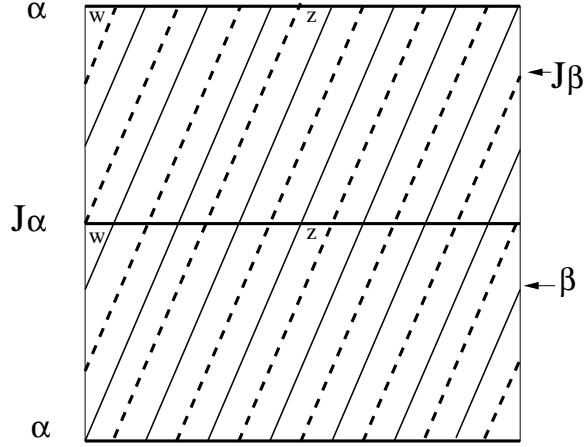, width=3in}
\end{center}
\caption{The twisted toroidal grid diagram which is a $4$-pointed genus $1$ Heegaard diagram for $\widetilde{K}_{7,3} \subset -L(7,3)$ (identify top-bottom, left-right, in the standard way).}
\label{fig:K7_3Simple}
\end{figure}

The two curves of slope $0$ on $T^2$ are the image in $T^2$ of the lines $y = 0$ and $y = \frac{1}{2}$ and the two curves of slope $\frac{p}{q}$ are the image in $T^2$ of the lines $y = \frac{p}{q}x$ and $y = \frac{p}{q}(x - \frac{1}{2})$.  We now identify the toroidal grid diagram with the fundamental domain $[0,1] \times [0,1] \subset \mathbb{R}^2$, \footnote{Note that by ``the image in $T^2$'' we mean the image of these lines in the quotient $\mathbb{R}^2/\mathbb{Z}^2$ and not the intersection of these lines with the chosen fundamental domain.}  and position our four basepoints at $$(\epsilon, 1-\epsilon), (\frac{1}{2} + \epsilon, 1-\epsilon), (\epsilon, \frac{1}{2} - \epsilon), (\frac{1}{2} + \epsilon, \frac{1}{2} - \epsilon)$$ where $0<\epsilon<\min(\frac{1}{p}, \frac{1}{q})$.  See Figure \ref{fig:K7_3Simple} for the example of $K_{7,3}$.

Following \cite{GT0607691} we construct a chain complex

\begin{itemize}
\item whose generators are indexed by bijections between the set of slope $0$ curves and the set of slope $\frac{p}{q}$ curves,
\item whose differentials are given by counting parallelograms,
\item and where the absolute Alexander (filtration) grading, absolute $\mathbb{Z}_2$ Maslov (homological) grading and Spin$^c$ structure of a generator is obtained by performing a sum of local gradings assigned to the vertices of the grid diagram.
\end{itemize}

Furthermore, by locating certain canonical generators and appealing to \cite{MR1957829}, we can improve the absolute $\mathbb{Z}_2$ to an absolute $\mathbb{Q}$ grading.

The paper is organized as follows.

In Section \ref{section:Background}, we fix notations and conventions and recall necessary Heegaard Floer homology background. 

In Section \ref{section:Diagram}, we describe the construction of a twisted toroidal grid diagram for the double-branched cover of a two-bridge knot.  We also show how to combinatorially read off Alexander, Spin$^c$, and $\mathbb{Q}$ Maslov gradings.

In Section \ref{section:Examples}, we perform a sample computation of the filtered chain complex in all Spin$^c$ structures for $\widetilde{K}_{3,1} \subset -L(3,1)$.

In Section \ref{section:Conjectures}, we briefly discuss generalizations of our methods.

{\bf Acknowledgements:} This paper was motivated by joint work in progress with Danny Ruberman and Sa{\v{s}}o Strle whose aim is to develop efficient calculational techniques with an eye towards obtaining new concordance invariants for knots.

I am indebted to Matt Hedden, Robert Lipshitz, Peter \Ozsvath, Danny Ruberman, Sa{\v{s}}o Strle and Jiajun Wang for many interesting conversations during the course of this project.  I am especially grateful to Peter \Ozsvath for patiently answering questions about \cite{GT0607691} and to Danny Ruberman for many helpful comments on a previous draft.

I would also like to remark that it was observed independently by Matt Hedden that general twisted toroidal grid diagrams (with $n$ parallel curves of slope $0$, $n$ parallel curves of slope $\frac{p}{q}$, and $2n$ basepoints) can be used to obtain a combinatorial description of the Heegaard Floer homology of any knot in a lens space.

\section{Heegaard Floer Homology Background}\label{section:Background}
We recall here the relevant constructions of Heegaard Floer homology groups associated to a nullhomologous knot $K$ in a closed, oriented, connected three-manifold $Y$.  Although many of these concepts work equally well for links, we will focus on the case of knots in this paper.  In addition, we will use $\mathbb{Z}_2$ coefficients in order to avoid talking about orientations of moduli spaces of holomorphic disks.  We also assume for technical reasons that $Y$ is a $\mathbb{Q}$ homology sphere.  We will, for simplicity's sake, focus on the $\widehat{HF}$ version of the theory, since this will be of the most interest to us for applications.  However, we remark that our method provides a combinatorial method for computing all versions of the theory in the case of $\widetilde{K}_{p,q} \subset \Sigma^2(K_{p,q})$.  For more details, see \cite{MR2113019, MR2065507, GT0512286, GT0306378}.

One specifies a knot $K$ in a three-manifold $Y$ by means of a $2n$-pointed Heegaard diagram (See Section 2 of \cite{GT0607691} as well as \cite{GT0512286}).  This data will be used to construct a filtered chain complex from which the filtered chain complex $\widehat{CF}(Y;K)$ (as well as the homology groups of the associate graded complex, $\widehat{HFK}(Y;K)$) can be extracted.

We begin by recalling how a $2n$-pointed Heegaard diagram specifies a knot in a three-manifold:

\begin{definition}
A $2n$-pointed Heegaard diagram compatible with an oriented knot $K$ in a three manifold $Y$ is a tuple $(S, \vec{\alpha}, \vec{\beta}, {\bf w}, {\bf z})$ associated to a handlebody decomposition of $Y$ coming from a generic self-indexing Morse function $f:Y \rightarrow \mathbb{R}$ with $|f^{-1}(0)| = |f^{-1}(3)| = n$.  After equipping $Y$ with a metric with respect to which $f$ is Morse-Smale, we obtain negative gradient flowlines of $f$ with respect to this metric.   Label the index $0$ critical points $a_1, \ldots a_n$ and the index $3$ critical points $b_1, \ldots b_n$.
\begin{itemize}
\item $S = f^{-1}(3/2)$ is a genus $g$ surface,
\item $\vec{\alpha} = (\alpha_1, \alpha_2, \ldots, \alpha_{g+(n-1)})$ is a collection of pairwise disjoint, embedded, closed curves on $S$.  These are the intersections of $S$ with flowlines from the index $3$ to index $1$ critical points of $f$,
\item Similarly, $\vec{\beta} = (\beta_1, \beta_2, \ldots, \beta_{g+(n-1)})$ are the intersections of $S$ with flowlines from the index $2$ to index $0$ critical points.
\item $\vec{w} = (w_1, \ldots, w_n)$ and $\vec{z} = (z_1, \ldots, z_n)$ are two $n$-tuples of points (all distinct) on $S - \vec{\alpha} - \vec{\beta}$, where $w_i$ specifies a unique flowline $\gamma_i$ from $b_i$ to $a_i$ (the one that intersects $S$ at $w_i$) and $z_i$ specifies a unique flowline $\eta_i$ from $b_{\sigma(i)}$ to $a_i$ ($\sigma$ some permutation of $\{1,\ldots, n\}$).
\end{itemize}

Then $K$ is uniquely determined by this data as the isotopy class of $${\bigcup}_{i=1}^n -\gamma_i \cup \eta_i$$

\end{definition}

\cite{GT0512286, GT0607691} associate to this data a chain complex $\widehat{CF}(S,\vec{\alpha},\vec{\beta},{\bf w},{\bf z})$ over $\mathbb{Z}_2[U_2, \ldots, U_n]$ freely generated by the intersection points between the subvarieties $$\mathbb{T}_{\vec{\alpha}} = \alpha_1 \times \ldots \times \alpha_{g+n-1}$$ and $$\mathbb{T}_{\vec{\beta}} = \beta_1 \times \ldots \times \beta_{g+n-1}$$ in $Sym^{g+n-1}(S)$ with differential given by $$\partial({\bf x}) = \Large{\sum}_{{\bf y} \in \mathbb{T}_{\alpha} \cap \mathbb{T}_{\beta}}\Large{\sum}_{\{\phi \in \pi_2({\bf x},{\bf y})|\mu(\phi) = 1, n_{w_1}(\phi) = 0\}} \#\left(\frac{\mathcal{M}(\phi)}{\mathbb{R}}\right) U_2^{n_{w_2}}\cdots U_n^{n_{w_n}} {\bf y}$$

Here we have used the notation in \cite{GT0607691}.

Let $\widehat{CF}(S,\vec{\alpha},\vec{\beta},{\bf w},{\bf z})$ denote the filtered chain complex associate to this data (the filtration is described in detail in Subsection \ref{subsection:Alexander}) and $\widehat{HF}(S,\vec{\alpha},\vec{\beta},{\bf w},{\bf z})$ its filtered chain homotopy type.

We will also make use of the chain complex with the same generators, but with coefficients in $\mathbb{Z}_2$ rather than $\mathbb{Z}_2[U_2, \ldots U_n]$ and a restricted differential.  Namely, define $\widehat{CFK}(S,\vec{\alpha},\vec{\beta},{\bf w},{\bf z})$ to be the chain complex with the same generators as above but with $$\partial_K({\bf x}) = \Large{\sum}_{{\bf y} \in \mathbb{T}_{\alpha} \cap \mathbb{T}_{\beta}}\Large{\sum}_{\{\phi \in \pi_2({\bf x},{\bf y})|\mu(\phi) = 1, n_{z_i} = n_{w_i}(\phi) = 0\}} \#\left(\frac{\mathcal{M}(\phi)}{\mathbb{R}}\right) {\bf y}.$$  We denote by $\widehat{HFK}(S,\vec{\alpha},\vec{\beta},{\bf w},{\bf z})$ the homology of this chain complex.

{\bf Remark:} When $n=1$, our $K$-compatible $2$-pointed Heegaard diagram for $Y$ yields $\widehat{CF}(Y;K)$, the filtered chain complex arising in the original formulation \cite{MR2065507, GT0306378} of knot Floer homology.  The homology of the associated graded complex of this filtered chain complex, denoted $\widehat{HFK}(Y;K)$, is called the {\it knot Floer homology} of $K$ in $Y$, and the $E^\infty$ term of the spectral sequence arising from the filtration is just the ordinary $\widehat{HF}(Y)$ obtained by forgetting the data of the knot.

Generators in $\widehat{CF}(S,\vec{\alpha},\vec{\beta},{\bf w},{\bf z})$ are assigned a bigrading $(n_1,n_2) \in \mathbb{Z}\times \mathbb{Q}$ whose first component is called the Alexander (filtration) grading and second component is called the Maslov (homological) grading.

\subsection{Maslov gradings}\label{subsection:Maslov}
The relative Maslov (homological) grading between two generators ${\bf x}, {\bf y}$ with non-empty $\pi_2({\bf x},{\bf y})$ is given by $$\mathcal{M}({\bf x}) - \mathcal{M}({\bf y}) = \mu(\phi) -\sum_i2n_{w_i}(\phi)$$ where $\phi \in \pi_2({\bf x},{\bf y})$, $$n_{w_i} = \#(\phi \cap V_{w_i})$$ is the algebraic intersection number with the subvariety  $V_{w_i} := \{w_i\} \times Sym^{g+n-2}(S)$, and $\mu(\phi)$ is the Maslov index.

This relative grading can be lifted to an absolute $\mathbb{Q}$ grading following \cite{MR1957829}.  In brief, one assigns to the unique generator of $\widehat{HF}(S^3)$ the $\mathbb{Q}$-grading $0$, and obtains the $\mathbb{Q}$ grading of any other generator by examining maps induced by cobordisms, a process which can be tricky in general.  

For the case of interest to us here ($\widetilde{K}_{p,q} \subset \Sigma^2(K_{p,q})$), however, we will be able to do this explicitly (see Subsection \ref{subsection:Maslov}) by using the inductive formula developed in Section 4.1 of \cite{MR1957829}.  This is the only part of the construction that is not strictly combinatorial.  However, we can combinatorially define a relative $\mathbb{Q}$-grading (using a formula developed by Lee and Lipshitz in \cite{GT0608001}) as well as an absolute $\mathbb{Z}_2$ Maslov grading.

Recall the following definitions:

\begin{definition} An {\it absolute $\mathbb{Z}_2$-valued homological grading} on $\widehat{CF}(Y)$, $Y$ a $\mathbb{Q}$ homology sphere, is an assignment, $$M_{\mathbb{Z}_2}: \mathbb{T}_{\alpha} \cap \mathbb{T}_{\beta} \rightarrow \mathbb{Z}_2$$subject to the following two conditions:

\begin{enumerate}
\item The relative $\mathbb{Z}_2$ grading is given by $$M_{\mathbb{Z}_2}({\bf x}) - M_{\mathbb{Z}_2}({\bf y}) = sgn({\bf x}) - sgn({\bf y})$$ where sgn({\bf x}) refers to the natural local orientation of $Sym^{g+(n-1)}(S)$ induced by ${\bf x} \in \mathbb{T}_{\alpha} \cap \mathbb{T}_{\beta}$.
\item The absolute $\mathbb{Z}_2$ is chosen so that $\chi(\widehat{HF}(Y)) = |H_1(Y)|.$  See Section 5.1 of \cite{MR2113020}.
\end{enumerate}
\end{definition}

\begin{definition} Given $\mathfrak{s}$ a Spin$^c$ structure on a $\mathbb{Q}$ homology sphere $Y$, the {\it correction term} $d_{\mathfrak{s}}(Y)$ is defined as the minimal $\mathbb{Q}$ grading of any non-torsion element in the image of $HF^{\infty}(Y,\mathfrak{s})$ in $HF^+(Y,\mathfrak{s})$.  See Section 4 of \cite{MR1957829}.

If $Y$ is an L-space (I.e., $Y$ has the Heegaard Floer homology of a lens space--see Definition 1.1 in \cite{MR2168576}), $d_{\mathfrak{s}}(Y)$ is the $\mathbb{Q}$-grading of the unique generator of $\widehat{HF}(Y;\mathfrak{s})$.
\end{definition}

\subsection{Alexander gradings}\label{subsection:Alexander}
$\widehat{CF}(S,\vec{\alpha},\vec{\beta},{\bf w}, {\bf z})$ is equipped with a filtration coming from an {\it Alexander grading} on the generators of the chain complex.  We will define this grading in two steps: first for the traditional knot Floer homology complex $\widehat{CF}(Y;K)$ and then for the knot complex associated to a $2n$-pointed Heegaard diagram for the knot.

The grading is obtained by using the evaluation of $c_1(\underline{\mathfrak{s}})$  on a capped-off Seifert surface for $K$ in $Y_0(K)$ ($0$-surgery on $K$) for $\underline{\mathfrak{s}} \in$ Spin$^c(Y_0(K))$, the set of {\it relative Spin$^c$ structures} for $K$ in $Y$.  Here we are using Turaev's identification (\cite{MR1484699}) of Spin$^c$ structures with homology classes of non-vanishing vector fields.

More specifically, Section 2.3 of \cite{MR2065507} and Section 2.6 of \cite{MR2113019}), describe how to split the chain complex $\widehat{CFK}(Y;K)$ associated to a traditional $2$-pointed Heegaard diagram into subcomplexes indexed by Spin$^c(Y_0(K))$.  Furthermore, they specify a map $$f: \mathbb{T}_{\vec{\alpha}} \cap \mathbb{T}_{\vec{\beta}} \rightarrow Spin^c(Y_0(K)),$$ assigning generators of $\widehat{CFK}(Y;K)$ to Spin$^c(Y_0(K))$ structures.

They go on to construct a splitting $$p_1 \times p_2: \xymatrix{Spin^c(Y_0(K)) \ar^{\cong}[r] & (Spin^c(Y) \times \mathbb{Z})},$$ where $p_1$ is given by restricting the vector field to $Y-K$ and taking the unique extension to $Y$, and $p_2$ is given by evaluating $c_1(\underline{\mathfrak{s}})$ on a capped off Seifert surface $\hat{F}$ for the knot: $$p_2(\underline{\mathfrak{s}}) = \frac{1}{2}<c_1(\underline{\mathfrak{s}}), [\hat{F}]>.$$  If $Y$ is a rational homology sphere, the map $p_2$ is independent of the choice of Seifert surface for the knot.

The absolute Alexander grading is then uniquely determined for a generator ${\bf x}$ in relative Spin$^c$ structure $\underline{\mathfrak{s}}_x$ by $${\bf A}({\bf x}) = \frac{1}{2}<c_1(\underline{\mathfrak{s}}_x),[\hat{F}]>.$$

{\flushleft {\bf Remark:} The relative Alexander grading of two generators in the same Spin$^c(Y)$ structure can be measured by looking at a disk connecting them.  I.e., given $${\bf x, y} \in \mathbb{T}_{\vec{\alpha}} \cap \mathbb{T}_{\vec{\beta}}$$ and $\phi \in \pi_2({\bf x},{\bf y})$ a homotopy class of disk connecting ${\bf x}$ to ${\bf y}$ in $Sym^{g}(S)$ we have $$\mathcal{A}({\bf x}) - \mathcal{A}({\bf y}) = n_z(\phi) - n_w(\phi).$$}

Bearing this in mind, we now proceed to define an Alexander grading on the chain complex associated to a $2n$-pointed Heegaard diagram for a knot with the property that if $\pi_2({\bf x}, {\bf y})$ is non-empty, then $$\mathcal{A}({\bf x}) - \mathcal{A}({\bf y}) = \left( \Large{\sum}_{i=1}^{n}n_{z_i}(\phi)\right) - \left(\Large{\sum}_{i=1}^{n}n_{w_i}(\phi) \right).$$

Furthermore, this will imply that $$\mathcal{A}({U_2^{a_2}\cdots U_n^{a_n} \bf x}) - \mathcal{A}({\bf x}) = \sum_{i=2}^n a_i.$$

Proposition 2.3 of \cite{GT0607691} explains how to do this. In brief, $\widehat{HFK}(S,\vec{\alpha},\vec{\beta},{\bf w},{\bf z})$ (recall that this is the homology of the $\widehat{CFK}$ chain complex with differential $\partial_K$) is related to $\widehat{HFK}(S^3;K)$ in a simple way.  The argument outlined in their proof works equally well for any nullhomologous knot in a rational homology sphere.

\begin{proposition} \label{prop:2npointed}
Let $(S, \vec{\alpha}, \vec{\beta}, {\bf w}, {\bf z})$ be a $2n$-pointed admissible Heegaard diagram compatible with a nullhomologous knot $K$ in a rational homology sphere $Y$.  Then $$\widehat{HFK}(S,\vec{\alpha}, \vec{\beta}, {\bf w}, {\bf z}) \cong \widehat{HFK}(Y;K) \otimes V^{\otimes (n-1)},$$ where $V$ is the vector space over $\mathbb{Z}_2$ with two generators, one in bigrading $(-1,-1)$ and the other in bigrading $(0,0)$.
\end{proposition}

The first coordinate of the bigrading denotes the filtration (Alexander) grading, and the second denotes the homological (Maslov) grading.

\begin{figure}
\begin{center}
\epsfig{file=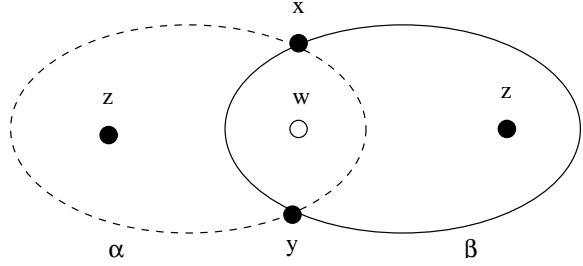, width=3in}
\end{center}
\caption{A $0-1$, $2-3$ Heegaard diagram stabilization for a knot $K$, obtained by replacing a small neighborhood of a point $p$ in $S- \vec{\alpha} - \vec{\beta} - \vec{w} - \vec{z}$ near $z$ with the above local picture.}
\label{fig:Newzandw}
\end{figure}

{\flushleft{\bf Proof of Proposition \ref{prop:2npointed}}
We refer the reader to the proof of Proposition 2.3 in \cite{GT0607691}, which works in this more general case.  The two main points are}
\begin{enumerate}
\item
The filtered chain homotopy type of an admissible Heegaard diagram with $2n$ basepoints (for a fixed $n$) is independent of the particular choice of diagram, since any two such diagrams can be connected by a sequence of isotopies, handleslides, and $1-2$ handle stabilizations/destabilizations.  Details about invariance under these moves can be found in \cite{MR2113019, MR2065507}.
\item
Adding an extra $w,z$ pair by introducing a pair of cancelling $0-1$ handle and $2-3$ handle pairs has the effect of tensoring the chain complex $\widehat{CFK}(S,\vec{\alpha},\vec{\beta},{\bf w},{\bf z})$ with $V$, since the introduction of such a pair is effected by replacing one of the $z$ basepoints with the local picture given in Figure \ref{fig:Newzandw}.  Notice that we get such a local picture by poking a trivial arc of the knot through the Heegaard surface at a point $p$ near $z$ (where by ``near z'' we simply mean that we can connect $z$ to $p$ by an arc on $S$ which does not intersect any $\alpha$ or $\beta$ curves).
\item Note that the new chain complex splits as the direct sum $C_x \oplus C_y$ of two copies of our old chain complex where $C_x$ consists of those generators of the form $(*,x)$ and $C_y$ consists of those generators of the form $(*,y)$.  Furthermore, two generators of this chain complex which agree everywhere except in the last component are connected by an obvious disk $\phi$ with $\sum n_{z_i}(\phi) - \sum n_{w_i}(\phi) = 1$ and $\mu(\phi) = 1$.  Iterating this process yields the desired conclusion.
\end{enumerate}

$\qed$

As mentioned earlier, the Alexander grading on $\widehat{CF}(S,\vec{\alpha},\vec{\beta},{\bf w},{\bf z})$ naturally gives rise to a $\mathbb{Z}$ filtration.  In particular, (see Section 3.1 of \cite{GT0607691}) one defines subcomplexes $\widehat{F}(K,m) \subset \widehat{CF}(S,\vec{\alpha},\vec{\beta},{\bf w},{\bf z})$ generated by elements $U_2^{a_2} \cdots U_n^{a_n} {\bf x}$ with $$\mathcal{A}(U_2^{a_2} \cdots U_n^{a_n}{\bf x}) = \mathcal{A}({\bf x}) - \sum_{i=2}^n a_i \leq m.$$

This filtration splits naturally into Spin$^c(Y)$ buckets.  Accordingly, one denotes by $\widehat{\mathcal{F}}_{\mathfrak{s}}(K,m)$ the subcomplex generated by elements in $\mathfrak{s}$ with $\mathcal{A}({\bf x}) \leq m$.

\begin{definition}
If $Y$ is an L-space, $\tau_{\mathfrak{s}}(K)$ is defined to be the minimal $m$ for which the map induced on homology $$i_*:H_*(\widehat{\mathcal{F}_{\mathfrak{s}}}(K,m)) \rightarrow \widehat{HF}(S,\vec{\alpha},\vec{\beta},{\bf w}, {\bf z})$$ is non-trivial.
\end{definition}

{\bf Remark:} $\tau_{\mathfrak{s}}(K)$ is defined in far greater generality (see Section 5 of \cite{MR2026543}), but we restrict to this special case for the sake of exposition, since it is all we will need at present.

\section{Nice Heegaard Diagrams for $(\Sigma^2(K_{p,q}), \widetilde{K}_{p,q})$}\label{section:Diagram}
In this section, we explicitly construct the $4$-pointed twisted toroidal grid diagram for $\widetilde{K}_{p,q} \subset \Sigma^2(K_{p,q})$ described in the introduction.

Once again, our convention is to denote by $K_{p,q}$ the two-bridge link with double-branched cover $-L(p,q)$.

\begin{figure}
\begin{center}
\epsfig{file=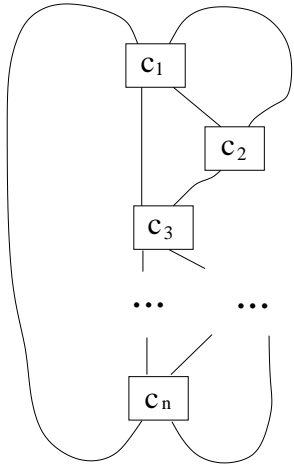, width=1.5in}
\end{center}
\caption{A two-bridge link with crossing numbers $(c_1, c_2, \ldots, c_n)$.  Positive = RH, and Negative = LH.}
\label{fig:Twobridge}
\end{figure}

Recall (see, e.g., Chapter 12 of \cite{MR1959408}) that a two-bridge link $L$ has the form given in Figure \ref{fig:Twobridge}.  It is a standard fact that two-bridge links are classified by the oriented homeomorphism class of their double-branched covers, the lens spaces $L(p,q)$ (where we assume that $p \in \mathbb{Z}_+, q \in \mathbb{Z}$, and $(p,q)=1$).  More precisely:

\begin{theorem} \cite{MR0082104, MR0258014,RTorsion,MR0116336} A two-bridge link with crossing numbers $$(c_1, c_2, \ldots, c_n)$$ has double-branched covering $-L(p,q)$, where $$\frac{p}{q} = c_1 - \frac{1}{\displaystyle c_2 -\frac{1}{\displaystyle c_3 - \frac{1}{\displaystyle \cdots - \frac{1}{c_n}}}}.$$  Furthermore, $K_{p,q}$ and $K_{p',q'}$ are isotopic iff $-L(p,q)$ is homeomorphic to $-L(p',q')$ by an orientation-preserving homeomorphism.  This condition is equivalent to the two conditions 

\begin{enumerate}
\item $p=p'$ and
\item $q = q' \mod p$ or $qq' = 1 \mod p$
\end{enumerate}

\end{theorem}

Note that $p$ is odd iff $L$ has one component.

We begin by describing the construction of a {\it nice} Heegaard knot diagram for the preimage of $K_{p,q}$ inside its double-branched cover.  An admissible Heegaard diagram (see \cite{MR2113019}) is said to be {\it nice} if all fundamental domains are either bigons or quadrilaterals.  In this case, the Floer homology is combinatorially defined.  See \cite{GT0607777, GT0607691}.

\begin{figure}
\begin{center}
\epsfig{file=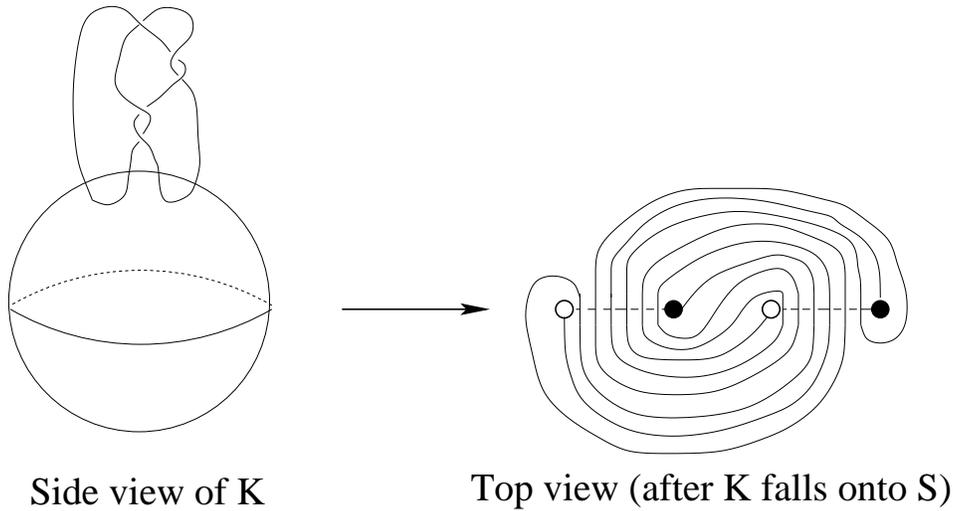, width=5in}
\end{center}
\caption{Obtaining the Schubert normal form of $K_{7,3}$}
\label{fig:SchubertSphere}
\end{figure}

Consider the Schubert normal form for the two-bridge knot $K_{p,q}$.  Recall that this is obtained by imbedding $K_{p,q}$ radially in $S^3$ so that $S$, the $S^2$ at radius $1$ from the origin, intersects $K_{p,q}$ in four points that are very close to the bottom two bridges.  Then let the knot fall onto $S$, keeping track of the over and under crossings.  When viewed on $S$, the knot is the union of two straight {\it underbridges} and two curvy {\it overbridges}.  See Figure \ref{fig:SchubertSphere}.

From the Schubert normal form for $K_{p,q}$ we construct a genus $0$, four-pointed Heegaard diagram for $S^3$ compatible with $K_{p,q}$. 

\begin{figure}
\begin{center}
\epsfig{file=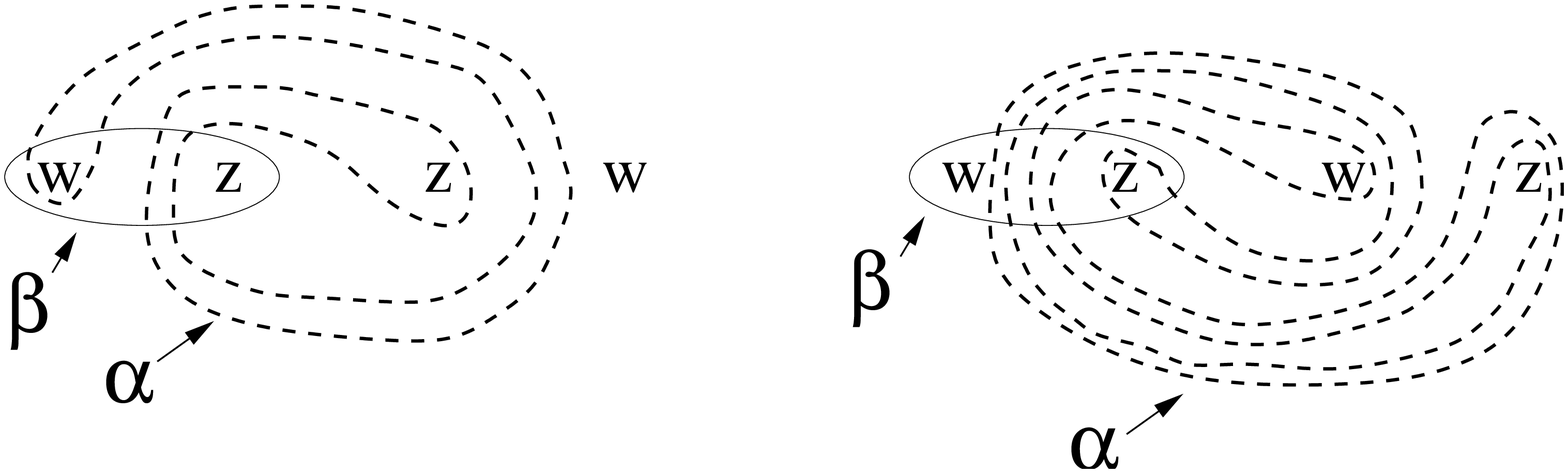, width=5in}
\end{center}
\caption{$4$-pointed Heegaard diagrams compatible with the knots $K_{3,1}$ and $K_{5,2}$, respectively}
\label{fig:EvenOdd}
\end{figure}

In the case where $K$ is the two-bridge knot $K_{p,q}$, we identify the maxima of the two upper bridges with the two index $0$ critical points of $f$ (labeled $a_1$ and $a_2$) and the minima of the two lower bridges with the two index $3$ critical points of $f$ (labeled $b_1$ and $b_2$).  Then the $\alpha$ curve is a regular neighborhood of either overbridge and the $\beta$ curve is a regular neighborhood of either underbridge in the Schubert normal form for $K_{p,q}$.  If $q$ is odd, then the points (read from left to right) are $w_1, z_1, z_2, w_2$, while if $q$ is even, then the points (read from left to right) are $w_1, z_1, w_2, z_2$.  See Figure \ref{fig:EvenOdd}.

\begin{figure}
\begin{center}
\epsfig{file=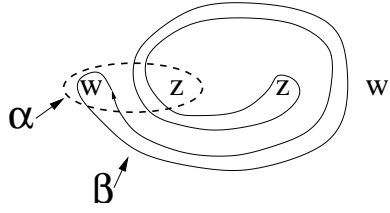, width=2in}
\end{center}
\caption{$4$-pointed Heegaard diagram compatible with $K_{3,1}$ after straightening out the $\alpha$ curve}
\label{fig:EvenOddAfter}
\end{figure}

Since it is traditional in the literature for the $\alpha$ curves to be simple and the $\beta$ curves complicated, we will straighten out the $\alpha$ curve by performing an isotopy of the $\alpha$ and $\beta$ curves avoiding the $w_i$ and $z_i$, which has the effect of reflecting the picture about a horizontal axis and swapping the roles of $\alpha$ and $\beta$.  Alternatively, we could have used the mirror, $K_{p,-q}$, and switched the roles of $\alpha$ and $\beta$ before letting the knot fall onto $S$ to get a Heegaard diagram compatible with $K_{p,q}$.  See Figure \ref{fig:EvenOddAfter}.

\begin{figure}
\begin{center}
\epsfig{file=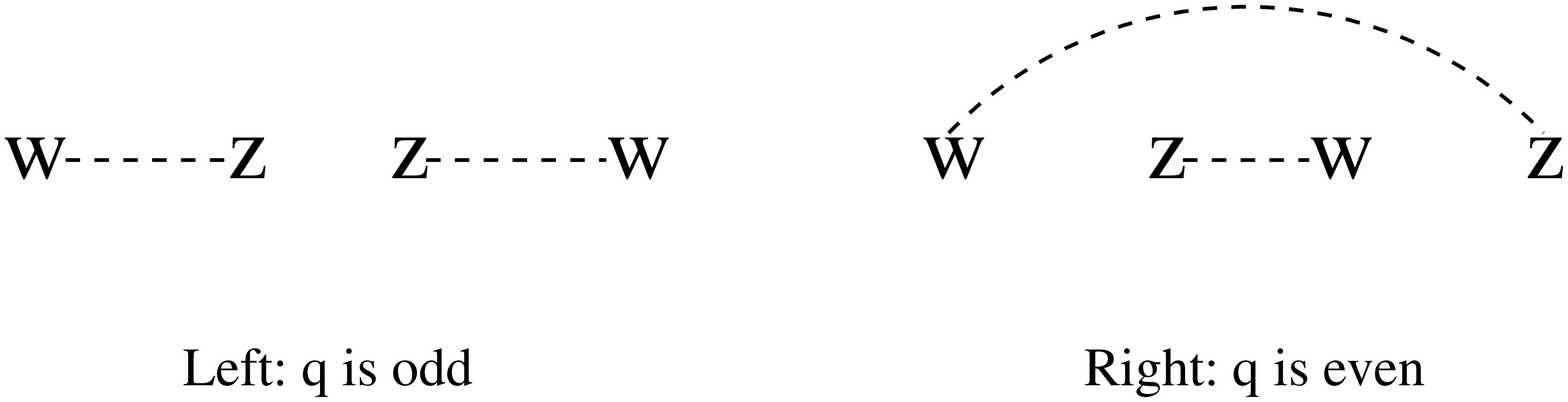, width=4in}
\end{center}
\caption{The intersection of a Seifert surface with $S$ (indicated by the dotted line) in the case where $q$ is odd (left) and even (right)}
\label{fig:SSurfaceIntersect}
\end{figure}

\begin{figure}
\begin{center}
\epsfig{file=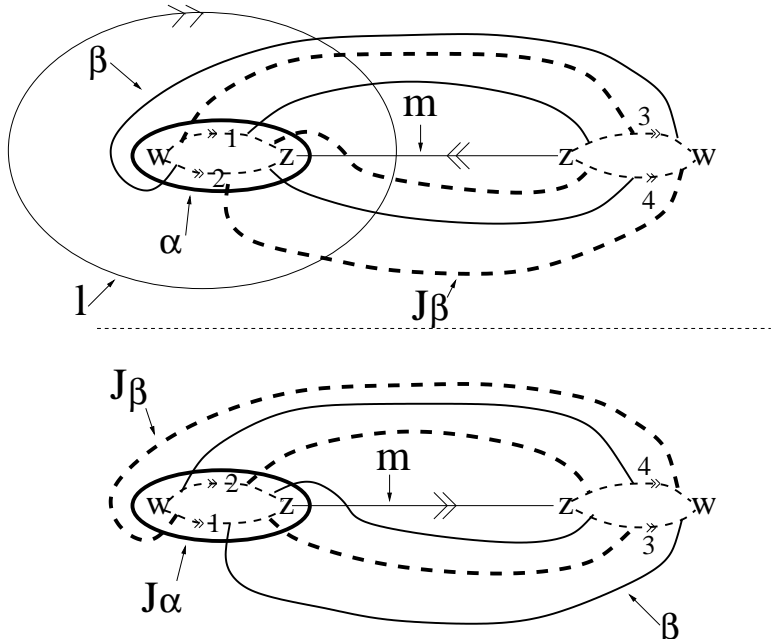, width=4in}
\end{center}
\caption{Form the Heegaard diagram, $\widetilde{S}$, for the double-branched cover, $\Sigma^2(K_{3,1}) = -L(3,1)$ compatible with $\widetilde{K}_{3,1}$ by gluing two copies of $S$ together along the branch cuts, yielding a torus.}
\label{fig:Schubert3_1}
\end{figure}

For the sake of simplicity, we will assume that $-p<q<p$ and $q$ is odd.  In this case, we form the Heegaard diagram for $\Sigma^2(K_{p,q})$ compatible with the preimage of $K_{p,q}$ by gluing two copies of this Heegaard diagram together along the two branch cuts obtained by examining where a Seifert surface for $K$ intersects $S$, as in Figure \ref{fig:Schubert3_1}.

If we take the curves ${\bf l}$ and ${\bf m}$ as a symplectic basis for $H_1(\widetilde{S})$, it is clear that $\alpha$ represents ${\bf l}$ and $\beta$ represents $q{\bf l} + p{\bf m}$.

\begin{figure}
\begin{center}
\epsfig{file=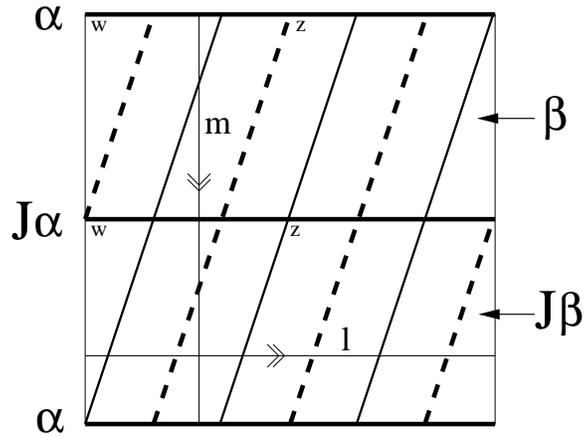, width=3in}
\end{center}
\caption{A square picture of the genus $1$ Heegaard diagram for $-L(3,1)$ compatible with $K_{3,1}$}
\label{fig:Square}
\end{figure}

A more convenient (square) picture of this diagram is given in Figure \ref{fig:Square}.

Notice that the fundamental domains in this diagram are all quadrilaterals, and, hence, all versions of the knot Floer homology are combinatorially defined.

We now proceed to characterize the generators and differentials of the chain complex, as well as the Alexander (filtration) gradings, Spin$^c$ structures, and absolute $\mathbb{Q}$ gradings of all of the generators.

\subsection{Generators and boundary maps}
The two $\alpha$ and two $\beta$ curves in $\widetilde{S}$ are the lifts of the single $\alpha$ and $\beta$ curve in the Heegaard diagram $S$ for $K_{p,q}$ in $S^3$.  We label the curves $\alpha$, $J\alpha$, $\beta$, and $J\beta$.  The $\mathbb{Z}_2$ symmetry coming from the double-branched cover construction gives rise to an involution, $J$, sending 

\begin{itemize}
\item $\alpha \leftrightarrow J(\alpha)$
\item $\beta \leftrightarrow J(\beta)$
\end{itemize}

The generators of the chain complex are intersection points, $\mathbb{T}_{\bf \alpha} \cap \mathbb{T}_{\bf \beta}$, in Sym$^2(\widetilde{S})$.  In other words, they are pairs of vertices in the twisted toroidal grid diagram associated to the graph of a bijection between the set of $\alpha$ and the set of $\beta$ curves.

The boundary map counts parallelograms of Maslov index one connecting ${\bf x}$ to ${\bf y}$.  As in \cite{GT0607691}, $\mu(\phi) = 1$ for $\phi \in \pi_2({\bf x},{\bf y})$ a parallelogram iff $$P_{\bf x}(\phi) + P_{\bf y}(\phi) = 1$$ where $P_{\bf x}(\phi)$ is the sum of the average local multiplicities of $\phi$ at the intersection points comprising ${\bf x}$.

\subsection{Alexander grading}
Note that the oriented knot itself can be seen in the diagram as a piecewise linear union of horizontal arcs and arcs of slope $\frac{p}{q}$ by connecting the basepoints, according to the conventions:

\begin{itemize}
\item The first horizontal arc is oriented from the $w$ to the $z$ on the top line,
\item Travel in a straight line until reaching another basepoint; then turn right.  Repeat the process until you return to where you started.
\end{itemize}

Note that $\widetilde{K}_{p,q}$ looks like a parallelogram when lifted to the universal cover, $\mathbb{R}^2$, of $\widetilde{S}$.  See Figure \ref{fig:K3_1Alexander}.

Again following \cite{GT0607691}, we define a local (relative) Alexander grading $${\bf A}: X \rightarrow \mathbb{Z}$$ on the set, $X$, of lattice points of our twisted toroidal grid diagram by counting the algebraic intersection number with $K$ of a straight line path $\gamma(x_1, y_1)$ from $x_1$ to $y_1$.  More succinctly,  $${\bf A}(y_1) - {\bf A}(x_1) = \#(K \cap \gamma(x_1,y_1)).$$

Figure \ref{fig:K3_1Alexander} gives such an assignment for the case of $\widetilde{K}_{3,1}$.

\begin{proposition} \label{proposition:AGradingEqual} The relative Alexander grading $\mathcal{A}({\bf x})$ for a chain complex generator ${\bf x} = (x_1, x_2)$ is obtained by summing over local gradings:  $$\mathcal{A}({\bf x}) = {\bf A}(x_1) + {\bf A}(x_2)$$
\end{proposition}

{\flushleft{\bf Proof of Proposition \ref{proposition:AGradingEqual}}  It is clear that if ${\bf x}, {\bf y}$ are two generators with $\pi_2({\bf x},{\bf y})$ non-empty (containing $\phi$), then the Alexander grading will match up with this local sum, for then $$\mathcal{A}({\bf x}) - \mathcal{A}({\bf y}) = \left(\sum_{i=1}^2 n_{z_i}(\phi)\right) - \left(\sum_{i=1}^2 n_{w_i}(\phi)\right),$$ and the local function ${\bf A}$ assigns to each element of $X$ the relative winding number of the lift of $K$ in $\mathbb{R}^2 = \widetilde{T}^2$ about that element.}

If $\pi_2({\bf x},{\bf y})$ is empty, then $\epsilon({\bf x},{\bf y})$ (which, we recall, is the element of $H_1(Y)$ represented by a path connecting ${\bf x}$ to ${\bf y}$ along $\alpha$ curves and ${\bf y}$ to ${\bf x}$ along $\beta$ curves) represents a non-trivial element of $H_1(Y)$.  

For each $a \in H_1(Y)$, choose a representative curve $\gamma_a$ on $\widetilde{S}$ with the property that $[\gamma_a] = a$ and $lk(\gamma_a, K) = 0$.  Then $\epsilon({\bf x},{\bf y}) - \gamma_{\epsilon({\bf x},\epsilon({y}))}$ represents $0$ in $H_1(Y)$, implying that $$\epsilon'({\bf x},{\bf y}) := \epsilon({\bf x},{\bf y}) - \gamma_{\epsilon({\bf x},{\bf y})} + n_{\alpha}\alpha + n_{\beta}\beta$$ represents $0$ in $H_1(T)$ for some $n_{\alpha}, n_{\beta} \in \mathbb{Z}$.

Notice that $$lk(\epsilon'({\bf x},{\bf y}), K) = lk(\epsilon({\bf x},{\bf y}), K)$$ since $\gamma_{\epsilon({\bf x},{\bf y})}, \alpha, \beta$ all have $0$ linking number with $K$.  But since $\epsilon'({\bf x},{\bf y})$ bounds a $2$-chain $\phi$ we can now measure $lk(\epsilon({\bf x},{\bf y}),K) = lk(\epsilon'({\bf x},{\bf y}),K)$ by counting $$\left(\sum_{i=1}^2 n_{z_i}(\phi)\right) - \left(\sum_{i=1}^2 n_{w_i}(\phi)\right).$$ Recalling that $$\epsilon({\bf x},{\bf y}) = \frac{1}{2}PD(c_1(\mathfrak{s}_x) - c_1(\mathfrak{s}_y))$$ completes the proof. $\qed$

One lifts to an absolute Alexander grading by making the unique choice yielding an appropriately symmetric Euler characteristic

$$\chi(\widehat{HFK}(Y;K,m)) := \sum_{d \in \mathbb{Z}_2} (-1)^d rk(\widehat{HFK}_d(Y;K,m))$$

where $$\widehat{HFK}(Y;K,m) := \bigoplus_{\{\underline{\mathfrak{s}}|m = \frac{1}{2}<c_1(\underline{\mathfrak{s}}),[\hat{F}]>\}} \widehat{HFK}(Y;K,\underline{\mathfrak{s}})$$ is the subcomplex of $\widehat{HFK}(Y;K)$ consisting of those Spin$^c$ structures with the appropriate evaluation on the homology class of a capped-off Seifert surface $\hat{F}$ for the knot.

More precisely, recall (see Proposition 3.10 of \cite{MR2065507}) that there is a conjugation symmetry on $\widehat{HFK}(Y;K)$ for $K$ a nullhomologous knot in a rational homology sphere.  This implies:

\begin{corollary}
$$\Delta_K(T) := \Large{\sum}_{i \in \mathbb{Z}} \chi(\widehat{HFK}(Y;K,m)) \cdot T^i$$ is a symmetric Laurent polynomial in $T$.
\end{corollary}

which allows us to fix the absolute Alexander grading by making the unique choice with the property $$\Large{\sum}_{{\bf x} \in \mathbb{T}_{\vec{\alpha}} \cap \mathbb{T}_{\vec{\beta}}} T^{{\bf A}({\bf x})} := \Delta_K(T) \cdot (1-T^{-1})^{n-1}.$$

\subsection{Spin$^c$ gradings}
We can make similar local assignments to partition the generators into Spin$^c(-L(p,q))$ structures.

More precisely, for the twisted grid diagram $\widetilde{S}$ associated to $\widetilde{K}_{p,q} \subset -L(p,q)$ we define a map $${\bf S}: X \rightarrow \mathbb{Z}_p$$ from intersection points to integers mod $p$.  Before describing how to do this, recall our identification in Section \ref{section:Introduction} of 

\begin{enumerate}
\item $\widetilde{S}$ with the fundamental domain $[0,1] \times [0,1]$,
\item $\alpha$ and $J(\alpha)$ with the images in $\widetilde{S} = \mathbb{R}^2/\mathbb{Z}^2$ of the lines $y=0(=1)$ and $y=\frac{1}{2}$, respectively,
\item $\beta$ and $J(\beta)$ with the images in $\widetilde{S}$ of the lines $y= \frac{p}{q}x$ and $y=\frac{p}{q}(x -\frac{1}{2})$ respectively.
\end{enumerate}

Now notice that the vertices of the twisted grid diagram are at $(\frac{j}{2p},0)$ and $(\frac{j}{2p}, \frac{1}{2})$ for $j \in 0, \ldots 2p-1$.  Furthermore, since (without loss of generality) we are restricting to the case that $q$ is odd and $-p<q<p$, we can identify the specific coordinates of the $4$ types of vertices $\alpha \cap \beta$, $\alpha \cap J(\beta)$, $J(\alpha) \cap \beta$, and $J(\alpha) \cap J(\beta)$.  Namely:

\vskip10pt

{\center\begin{tabular}{|r|c|c|}
\hline
$\cap $ & $\beta$ & $J(\beta)$ \\
\hline
$\alpha$ & $(\frac{j}{p},0)$ & $(\frac{1}{2p} + \frac{j}{p},0)$\\
\hline
$J(\alpha)$ & $(\frac{1}{2p} + \frac{j}{p}, \frac{1}{2})$ & $(\frac{j}{p}, \frac{1}{2})$\\
\hline
\end{tabular}}

\vskip10pt
{\flushleft for $j = 0, 1, \ldots p-1$.}

In the case $q>0$ we make the assignment:

\begin{enumerate}
\item ${\bf S}(\frac{j}{p},0) = j \in \mathbb{Z}_p$
\item ${\bf S}(\frac{1}{2p}+\frac{j}{p},0) = j \in \mathbb{Z}_p$
\item ${\bf S}(\frac{j}{p},\frac{1}{2}) = j + \frac{q-1}{2} \in \mathbb{Z}_p$
\item ${\bf S}(\frac{1}{2p} + \frac{j}{p},\frac{1}{2}) = j + \frac{q+1}{2} \in \mathbb{Z}_p$
\end{enumerate}

whereas in the case $q<0$ we make the assignment:

\begin{enumerate}
\item ${\bf S}(\frac{j}{p},0) = j \in \mathbb{Z}_p$
\item ${\bf S}(\frac{1}{2p}+\frac{j}{p},0) = j + 1 \in \mathbb{Z}_p$
\item ${\bf S}(\frac{j}{p},\frac{1}{2}) = j + \frac{q+1}{2} \in \mathbb{Z}_p$
\item ${\bf S}(\frac{1}{2p} + \frac{j}{p},\frac{1}{2}) = j +\frac{q+1}{2} \in \mathbb{Z}_p$
\end{enumerate}

See Figure \ref{fig:K3_1Spinc} in Section \ref{section:Examples} for the example of $K_{3,1}$.

\begin{proposition}\label{proposition:SpincMatch}
Let ${\bf x}$ be a generator of $\widehat{CF}(T)$ whose components are the intersection points $x_1$ and $x_2$. The function $$\mathcal{S}({\bf x}) := \sum_{i=1}^2 {\bf S}(x_i) \in \mathbb{Z}_p$$ partitions the generators of $\widehat{CF}(S,\vec{\alpha},\vec{\beta},{\bf w},{\bf z})$ into Spin$^c(Y)$ structures, $\mathfrak{s}_i$, indexed by elements of $\mathbb{Z}_p$.  

Furthermore, if $J({\bf x})$ is the image of ${\bf x}$ under the involution on $\widetilde{S}$ sending $\alpha \leftrightarrow J(\alpha)$ and $\beta \leftrightarrow J(\beta)$ and ${\bf x} \in \mathfrak{s}_{i}$, then $$J({\bf x}) \in \mathfrak{s}_{-i}.$$  The $J$-invariant generators are precisely the ones in $\mathfrak{s}_0$.
\end{proposition}

{\flushleft {\bf Proof of Proposition \ref{proposition:SpincMatch}}
First, recall that Spin$^c(-L(p,q))$ is an affine set for the action of $$H_1(-L(p,q)) \cong \mathbb{Z}_p.$$  Furthermore, since $p$ is odd, there is a unique spin Spin$^c$ structure, and it is natural to identify this Spin$^c$ structure with the $0$ element of $\mathbb{Z}_p$.  This choice allows us, once we have chosen a generator $\gamma \in H_1(-L(p,q))$, to identify the Spin$^c$ structures with $\mathbb{Z}_p$ by comparing $c_1(\mathfrak{s})$ with $\gamma$.}

More precisely, if $PD(c_1(\mathfrak{s})) = 2m \cdot \gamma$, then we identify $\mathfrak{s}$ with $m \in \mathbb{Z}_p$.

To see that the given function partitions generators of $\widehat{CF}(T)$ in the desired fashion, begin by noticing that the $J$-invariant generators, ${\bf x}= (x,J(x))$, are precisely the ones for which $$\sum_{i=1}^2{\bf S}(x_i) = 0 \in \mathbb{Z}_p,$$ since $J$ can be described as the $180$ degree rotation in $\mathbb{R}^2$ about the center of the cell containing the basepoint $w$ in the upper left-hand corner \footnote{Note that there is nothing special about this basepoint.  A $180$-degree rotation about any of the four basepoints yields the same result.} and we have chosen our assignments ${\bf S}$ to be anti-symmetric with respect to this operation.  I.e., $${\bf S}(x_i) \equiv -{\bf S}(Jx_i) \mod p.$$  This proves that $$(x_1,x_2) \mbox{ is $J$-invariant} \Longrightarrow \sum_{i=1}^2 {\bf S}(x_i) = 0.$$

To measure the Spin$^c$ structure of any non $J$-invariant generator ${\bf y} = (y_1,y_2)$, first choose a curve $\gamma$ on $T$, representing a homology class on $T$ with the property $$<\gamma,\beta> = 1$$ to be our generator of $H_1(-L(p,q))$.  By $<-,->$, we mean the standard intersection pairing on $H_1(T)$.

Then simply connect ${\bf y}$ to ${\bf y}_0 := (y_1,Jy_1) \in \mathfrak{s}_0$ by paths along the $\alpha$ and $\beta$ curves to form $\epsilon({\bf y},{\bf y}_0).$  Notice that our assignment ${\bf S}$ was made so that $\epsilon({\bf y},{\bf y}_0)$ measures the intersection pairing with $\beta$: $$\mathcal{S}({\bf y}) = <\epsilon({\bf y},{\bf y}_0),\beta>,$$ yielding our chosen isomorphism of $H_1(-L(p,q))$ with $\mathbb{Z}_p$. 

$\qed$

\subsection{Absolute $\mathbb{Z}_2$ and $\mathbb{Q}$ homological gradings}
We can specify the absolute $\mathbb{Z}_2$ homological grading on generators by performing a similar local sum.

Specifically, recall that we are assuming that $-p<q<p$, $q$ is odd, and we are specifying the intersection points in terms of their coordinates in the fundamental domain $[0,1] \times [0,1] \subset \mathbb{R}^2$.  Then define $${\bf M}: X \rightarrow \left\{0, \frac{1}{2}\right\}$$ by the simple rule:

If $q>0$,

\begin{enumerate}
\item ${\bf M}(\frac{i}{p},0) = {\bf M}(\frac{i}{p}, \frac{1}{2}) = 0$
\item ${\bf M}(\frac{1}{2p} + \frac{i}{p},0) = {\bf M}(\frac{1}{2p} + \frac{i}{p},\frac{1}{2}) = \frac{1}{2}$ 
\end{enumerate}

and if $q<0$,

\begin{enumerate}
\item ${\bf M}(\frac{i}{p},0) = {\bf M}(\frac{i}{p},\frac{1}{2}) = \frac{1}{2}$
\item ${\bf M}(\frac{1}{2p} + \frac{i}{p},0) = {\bf M}(\frac{1}{2p} + \frac{i}{p},\frac{1}{2}) = 0$ 
\end{enumerate}

\begin{proposition}\label{prop:Z2grading}
We can calculate the absolute $\mathcal{M}_{\mathbb{Z}_2}$ grading of a generator by summing over these local gradings.  I.e., if ${\bf x}$ consists of the intersection points $x_1$ and $x_2$ then $$\mathcal{M}_{\mathbb{Z}_2}({\bf x}) = \sum_{i=1}^2 {\bf M}(x_i) \in \mathbb{Z}_2.$$
\end{proposition}

{\flushleft {\bf Proof of Proposition \ref{prop:Z2grading}:} This gives the correct relative $\mathbb{Z}_2$ grading because it gives the correct relative local intersection parity of the generators in $\mathbb{T}_{\vec{\alpha}} \cap \mathbb{T}_{\vec{\beta}}$ (see \cite{MR1241874}). Namely, notice that our assignment insures that the generators of the form $(\alpha \cap \beta, J\alpha \cap J\beta)$ are all in different $\mathbb{Z}_2$ Maslov grading than those of the form $(J\alpha \cap \beta, \alpha \cap J\beta)$, which agrees with the observation that the $J$ action flips local orientation.}

\begin{figure}
\begin{center}
\epsfig{file=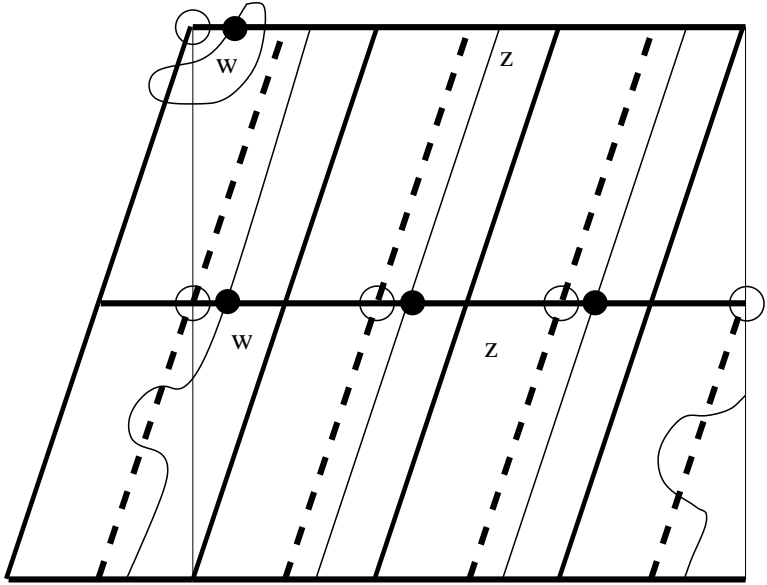, width=4in}
\end{center}
\caption{Finding homological grading $0$ generators in each Spin$^c$ structure}
\label{fig:Cangen}
\end{figure}

To see that this assignment gives the correct lift to an absolute $\mathbb{Z}_2$ grading, notice that we can locate a generator lying in $0$ homological grading in each Spin$^c$ structure by ``forgetting the knot.''  More precisely, by sliding $\beta$ over $J\beta$ and allowing ourselves to perform isotopies that cross the $z$ basepoints, we arrive at the diagram given in Figure \ref{fig:Cangen}.  Note that for each Spin$^c$ structure there is a Maslov index $0$ triangle joining the standard top-most generators of $\widehat{HF}(-L(p,q)) \otimes V$, labeled by the smaller black dots, to the generators in our original Heegaard diagram labeled by the larger white dots.  This agrees with our assignment of these generators to absolute $\mathbb{Z}_2$ grading $0$.  Note that one constructs a similar diagram in the case $q<0$, arriving at the opposite local assignment.  $\qed$

We are unable to specify an absolute $\mathbb{Q}$-grading by means of this type of local sum.  However, we can appeal to \cite{MR1957829}, in which \Ozsvath and \Szabo inductively compute the $\mathbb{Q}$-grading correction terms $d_{\mathfrak{s}}(-L(p,q))$ for lens spaces.  

The Maslov index zero triangles mentioned in the proof above connect the standard lens space generators (whose $\mathbb{Q}$ gradings are calculated in \cite{MR1957829}) to generators on our Heegaard diagram, yielding a $\mathbb{Q}$ grading assignment for one generator in each Spin$^c$ structure.  

Then we note, as in \cite{GT0607691}, that if ${\bf x},{\bf y}$ are generators in the same Spin$^c(Y)$ structure and $\phi \in \pi_2({\bf x},{\bf y})$ then $$\mathcal{M}({\bf x}) - \mathcal{M}({\bf y}) = P_{\bf x}(\phi) + P_{\bf y}(\phi)-2 \cdot W(\phi),$$ where $$P_{\bf x}(\phi) := \sum_{i=1}^2 p(x_i)(\phi),$$ $p_{x_i}(\phi)$ is the average of the local multiplicities of $\phi$ in the four quadrants around $x_i$ (See the definition of $n_p(A)$ in Section 4.2 of \cite{SG0502404}), and $$W(\phi) := \sum_{i=1}^2 n_{w_i}(\phi).$$

This gives us a relative $\mathbb{Z}$ grading on the generators in a single Spin$^c(Y)$ structure, yielding absolute $\mathbb{Q}$-gradings for all generators. 

\section{Example: $\widetilde{K}_{3,1} \subset -L(3,1)$}\label{section:Examples}

\begin{figure}
\begin{center}
\epsfig{file=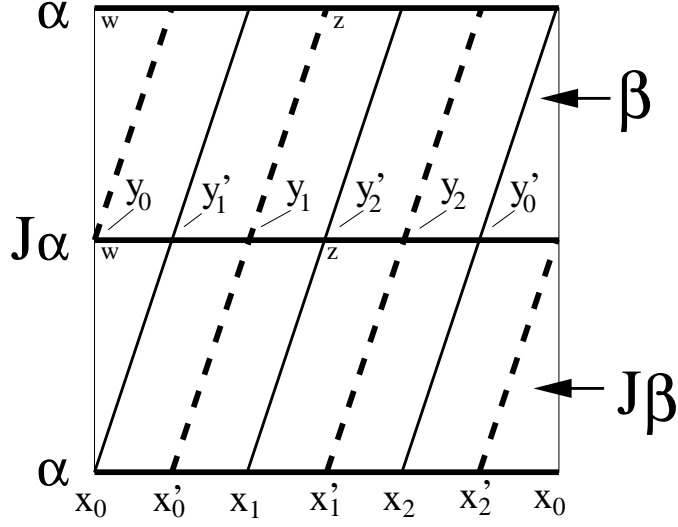, width=3.5in}
\end{center}
\caption{Genus $1$ Heegaard diagram for $\widetilde{K}_{3,1} \subset -L(3,1)$.  The intersection points have been labeled according to local Spin$^c$ grading, ${\bf S}$.}
\label{fig:K3_1Spinc}
\end{figure}

The genus $1$ Heegaard diagram for $\widetilde{K}_{3,1} \subset -L(3,1)$ is shown in Figure \ref{fig:K3_1Spinc}  Notice that there are $3$ intersections between each pair of $\alpha$, $\beta$ curves, yielding $3^2 + 3^2 = 18$ total generators in the chain complex $\widehat{CF}(S,(\alpha,J\alpha),(\beta,J\beta),{\bf w},{\bf z})$.  For convenience, we have labeled the intersection points according to their local Spin$^c$ grading.  Namely, 

\begin{enumerate}
\item $(\frac{i}{3},0) \in \alpha \cap \beta$ is labeled $x_i$
\item $(\frac{1}{6} + \frac{i}{3},0) \in \alpha \cap J\beta$ is labeled $x_i'$\item $(\frac{i}{3},\frac{1}{2}) \in J\alpha \cap J\beta$ is labeled $y_{i+1}$
\item $(\frac{1}{6}+ \frac{i}{3},\frac{1}{2}) \in J\alpha \cap \beta$ is labeled $y_{i + 2}'$.
\end{enumerate}

We remark that we have also chosen the labeling on the intersection points so that $J(x_i) = y_{-i \mod p}$ and $J(x_i') = y'_{-i \mod p}$.

\begin{figure}
\begin{center}
\epsfig{file=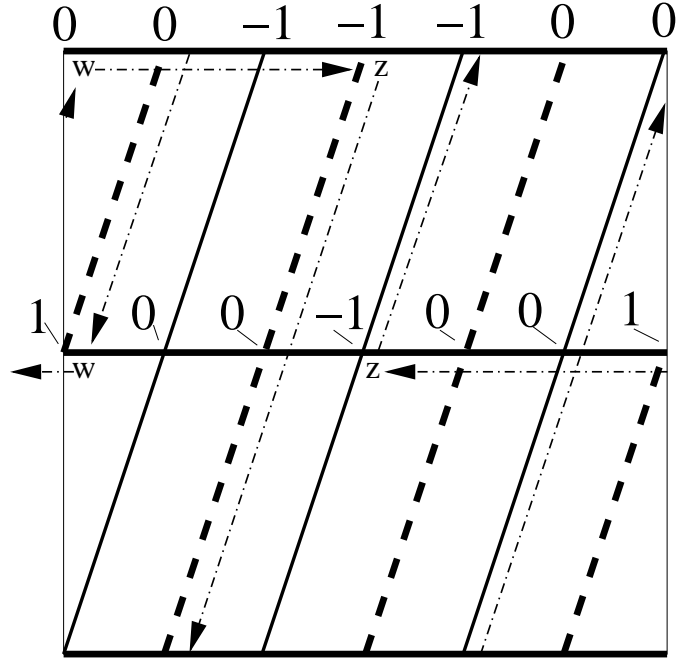, width=3.5in}
\end{center}
\caption{Here, the knot (dotted line) is allowed to fall on the genus $1$ Heegaard diagram for $\widetilde{K}_{3,1} \subset -L(3,1)$.  The number at each intersection point corresponds to the local Alexander grading, ${\bf A}$.}
\label{fig:K3_1Alexander}
\end{figure}

In Figure \ref{fig:K3_1Alexander} we have made assignments for the local Alexander grading ${\bf A}$.  The absolute Alexander grading of a generator is obtained by summing the local contribution at each intersection point.  So, for example, the intersection point $(x_0,y_1)$ has absolute Alexander grading $1$.

\begin{figure}
\begin{center}
\epsfig{file=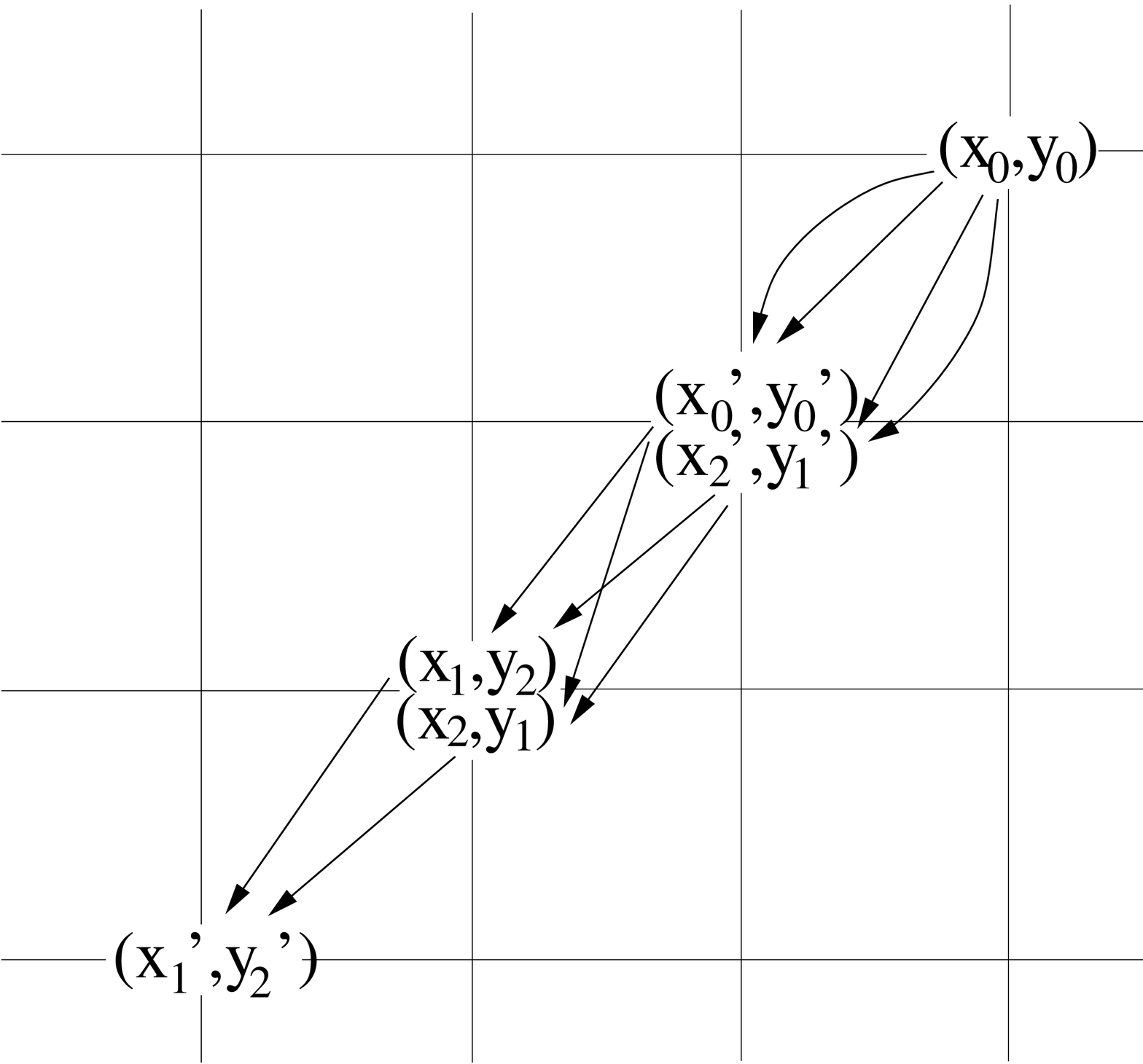, width=3in}
\end{center}
\caption{Generators and differentials in $\widehat{CF}(S,(\alpha,J\alpha),(\beta,J\beta),{\bf w},{\bf z};\mathfrak{s}_0)$.}
\label{fig:K3_1_CF_s0}
\end{figure}

\begin{figure}
\begin{center}
\epsfig{file=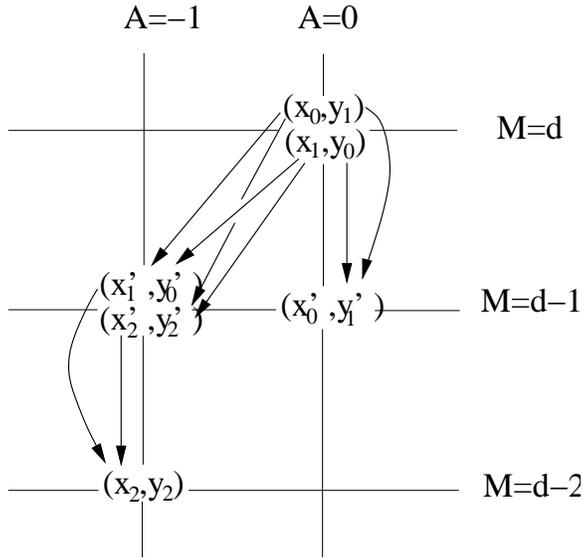, width=3in}
\end{center}
\caption{Generators and differentials in $\widehat{CF}(S,(\alpha,J\alpha),(\beta,J\beta),{\bf w},{\bf z};\mathfrak{s}_1)$.}
\label{fig:K3_1_CF_s1}
\end{figure}

We proceed, in Figures \ref{fig:K3_1_CF_s0} and \ref{fig:K3_1_CF_s1}, to calculate the filtered chain homotopy type of the complex $\widehat{CF}(S,\vec{\alpha},\vec{\beta},{\bf w},{\bf z})$ in all three Spin$^c$ structures.  Note that the chain complex for $\mathfrak{s}_{-1}$ is the same as the one for $\mathfrak{s}_1$, so we need only perform the calculation for $\mathfrak{s}_0$ and $\mathfrak{s}_1$.

Recall that the generators in $\mathfrak{s}_k$ are those $(x_i,y_j)$ and $(x_i',y_j')$ satisfying $i+j \equiv k \mod 3$.

In an attempt to keep the diagrams uncluttered, we have only included one copy of each generator in its appropriate bigrading (and not its $U_2$ translates).  Therefore, the homology (forgetting the filtration) will be $V = (\mathbb{Z}_2)_{(0,0)} \oplus (\mathbb{Z}_2)_{(-1,-1)}$ where the indices denote the Alexander-Maslov bigrading.

We easily compute the homology in $\mathfrak{s}_0$ as:

\begin{eqnarray*}
\widehat{HFK}(S,(\alpha,J\alpha),(\beta,J\beta),{\bf w},{\bf z};\mathfrak{s}_0) &=& (\mathbb{Z}_2)_{(-2,d_{}-3)} \oplus (\mathbb{Z}_2^2)_{(-1,d_{}-2)} \oplus (\mathbb{Z}_2^2)_{(0,d_{}-1)} \oplus (\mathbb{Z}_2)_{(1,d_{})}\\
&=& [(\mathbb{Z}_2)_{(-1,d_{}-2)} \oplus (\mathbb{Z}_2)_{(0,d_{}-1)} \oplus (\mathbb{Z}_2)_{(1,d_{})}] \otimes V 
\end{eqnarray*}

where $d = -\frac{1}{2}$ is the $\mathbb{Q}$ grading correction term for $\mathfrak{s}_0$.  In other words, $$\widehat{HFK}(-L(3,1);\widetilde{K}_{3,1};\mathfrak{s}_0) \cong \widehat{HFK}(S^3;K_{3,1})$$ (modulo an overall $\mathbb{Q}$ grading shift).  Note that $\tau_{\mathfrak{s}_0} = 1$.  Compare \cite{GT0507498}.

In $\mathfrak{s}_{\pm 1}$ we get:

\begin{eqnarray*}
\widehat{HFK}(S,(\alpha,J\alpha),(\beta,J\beta),{\bf w},{\bf z};s_{\pm 1}) &=& (\mathbb{Z}_2)_{(-1,d-1)} \oplus (\mathbb{Z}_2)_{(0,d)}\\
&=& (\mathbb{Z}_2)_{(0,d)} \otimes V 
\end{eqnarray*}

where $d= \frac{1}{6}$ is the correction term for $\mathfrak{s}_{\pm 1}$.  So $$\widehat{HFK}(-L(3,1);\widetilde{K}_{3,1};\mathfrak{s}_{\pm 1}) \cong \widehat{HFK}(S^3;U)$$ (modulo an overall $\mathbb{Q}$ grading shift).  Note that $\tau_{\mathfrak{s}_{\pm 1}} = 0$.

\section{Related Constructions} \label{section:Conjectures}
It is natural to ask how far the methods discussed here will allow us to go in calculating Floer homology invariants for higher cyclic branched covers of more general knots in $S^3$.

The first thing one might try is to take a general knot in $n$-bridge position, intersect with a transverse $S^2$ and once again form the $2n$-pointed genus $0$ Heegaard diagram for $K$ in $S^3$ obtained by letting the strands of the knot fall down on $S^2$ and using a collection of regular neighborhoods of the strands as the $\alpha$ and $\beta$ curves.

Unfortunately, this will not yield a ``nice'' Heegaard diagram when $n>2$ (there will be fundamental domains which are not quadrilaterals or bigons).  A suitable application of the methods of Sarkar and Wang \cite{GT0607777}, however, may yield good results.

In another direction, one can always combinatorially compute the most basic version of the knot Floer homology in any $m$-cyclic branched cover, i.e., the associated graded complex $\widehat{HFK}(\Sigma^m(K);\widetilde{K})$, by taking the $n$-grid presentation of $K$ in $S^3$ described in \cite{GT0607691} and branching around the $2n$ points.

For the $m$-fold cyclic cover this yields a $2n$-pointed genus $(m-1)n + 1$ Heegaard diagram for $\widetilde{K} \subset \Sigma^m(K).$ with $mn$ $\alpha$ curves and $mn$ $\beta$ curves, yielding $mn!$ generators for the Heegaard Floer chain complex.  The symmetries in this chain complex may be exploitable to make its homology nearly as fast to compute as the original $S^3$ knot Floer homology.

We should remark that in this case the fundamental domains containing the basepoints will all be $2mn$-gons, so we will only have combinatorial access to information about the associated graded complex, $\widehat{HFK}$.
\bibliography{TwoBridge_Tau}
\end{document}